\documentclass[final,letterpaper]{siamart0516}
\usepackage{adjustbox}
\usepackage{enumitem}
\usepackage{amssymb}
\usepackage{amsmath}
\usepackage{tabularx,multirow}
\usepackage{picture,slashbox}
\usepackage{graphicx}
\usepackage{epsfig,epsf,psfrag,epstopdf,subfigure}
\usepackage{color,bm}
\usepackage{comment}
\usepackage{algorithm,algorithmic}
\usepackage{url}
\usepackage{booktabs}
\usepackage{tablefootnote}

\newenvironment{myproof}{\paragraph{Proof}}{\hfill$\square$}

\newcommand{\be}{\begin{equation}}
\newcommand{\ee}{\end{equation}}
\newcommand{\ba}{\begin{array}}
\newcommand{\ea}{\end{array}}
\newcommand{\bea}{\begin{eqnarray}}
\newcommand{\eea}{\end{eqnarray}}
\newcommand{\beas}{\begin{eqnarray*}}
\newcommand{\eeas}{\end{eqnarray*}}

\newcommand{\bx}{{\textbf x}}
\newcommand{\by}{{\textbf y}}
\newcommand{\bz}{{\textbf z}}
\newcommand{\bn}{{\textbf n}}

\newcommand{\bh}{{\textbf h}}
\newcommand{\bp}{{\textbf p}}

\newcommand{\bk}{{\textbf k}}
\newcommand\bR{\textbf{R}}

\renewcommand{\theequation}{\arabic{section}.\arabic{equation}} %
\newtheorem{exmp}{Example}
\newtheorem{remark}{Remark}[section]

\newtheorem{thm}{Theorem}

\title{Fast convolution solver based on far-field smooth approximation \thanks{This
work was supported by the National Natural Science Foundation of China (No.12271400).}}

\author{Xin Liu\thanks{Center for Applied Mathematics and KL-AAGDM, Tianjin University, Tianjin 300072, China. ({\tt liuxin\_921@tju.edu.cn} (X. Liu); {\tt Zhang\_Yong@tju.edu.cn} (Y. Zhang).)} \and
Yong Zhang\footnotemark[2] \thanks{Corresponding author.} }

\begin{document}
\maketitle

\begin{abstract}
The convolution potential arises in a wide variety of application areas, and its efficient and accurate evaluation encounters
three challenges: singularity, nonlocality and anisotropy.
We introduce a fast algorithm based on a far-field smooth approximation of the kernel,
where the bounded domain Fourier transform, one of the most essential difficulties, is well approximated by the whole space Fourier transform which usually admits explicit formula.
The convolution is split into a regular and singular integral, and they are well resolved by trapezoidal rule and Fourier spectral method respectively.
 The scheme is simplified to a discrete convolution and is implemented efficiently with Fast Fourier Transform (FFT).
Importantly, the tensor generation procedure is quite simple, highly efficient and independent of the anisotropy strength.
It is easy to implement and achieves spectral accuracy with nearly optimal efficiency and minimum memory requirement.
Rigorous error estimates and extensive numerical investigations, together with a comprehensive comparison, showcase its superiorities for different kernels.
\end{abstract}


\pagestyle{myheadings}\thispagestyle{plain}
\markboth{X. Liu  and Y. Zhang}
{FSA-based fast convolution}


\section{Introduction}
Nonlocal potentials, which are given by a convolution of a translational invariant Green function with a
fast-decaying smooth function, are common and have wide-ranging applications.
Examples include the Newtonian potential in cosmology, the Poisson potential in electrostatics, plasma physics and quantum physics \cite{Bao2013Mathematical}. The efficient and accurate calculation of these nonlocal potentials is a prominent and vital area of research in the science and engineering community.
In this paper, we focus on the evaluation of convolution-type nonlocal potential
\be\label{convolution}
\Phi(\bx) =\int_{\mathbb{R}^d} U(\bx-\by) \rho(\by) {\rm d} \by
=\frac{1}{(2 \pi)^d} \int_{\mathbb{R}^d} \widehat{U}(\bk)  \widehat{\rho}(\bk) e^{i \bk \cdot \bx}{\rm d}\bk
,\quad  \bx \in \mathbb{R}^d,
\ee
where $d$ is the ambient dimension, the density $\rho(\bx)$ is a fast-decaying smooth function,
the kernel $U(\bx)$, a given radially symmetric function, is usually singular at the origin,
and $\widehat{f}(\bk)= \int_{\mathbb{R}^d} f(\bx) e^{-i \bk \cdot \bx} {\rm d} \bx$ is Fourier transform of $f(\bx)$.

\begin{figure}[!htbp]
\label{anisoDomain}
\centering
\includegraphics[width=5.1cm,height=1.4cm]{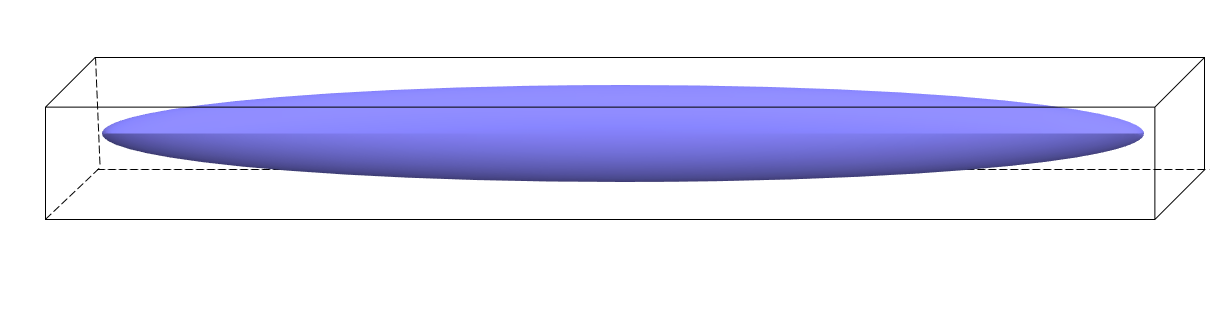}
\includegraphics[width=5.3cm,height=1.8cm]{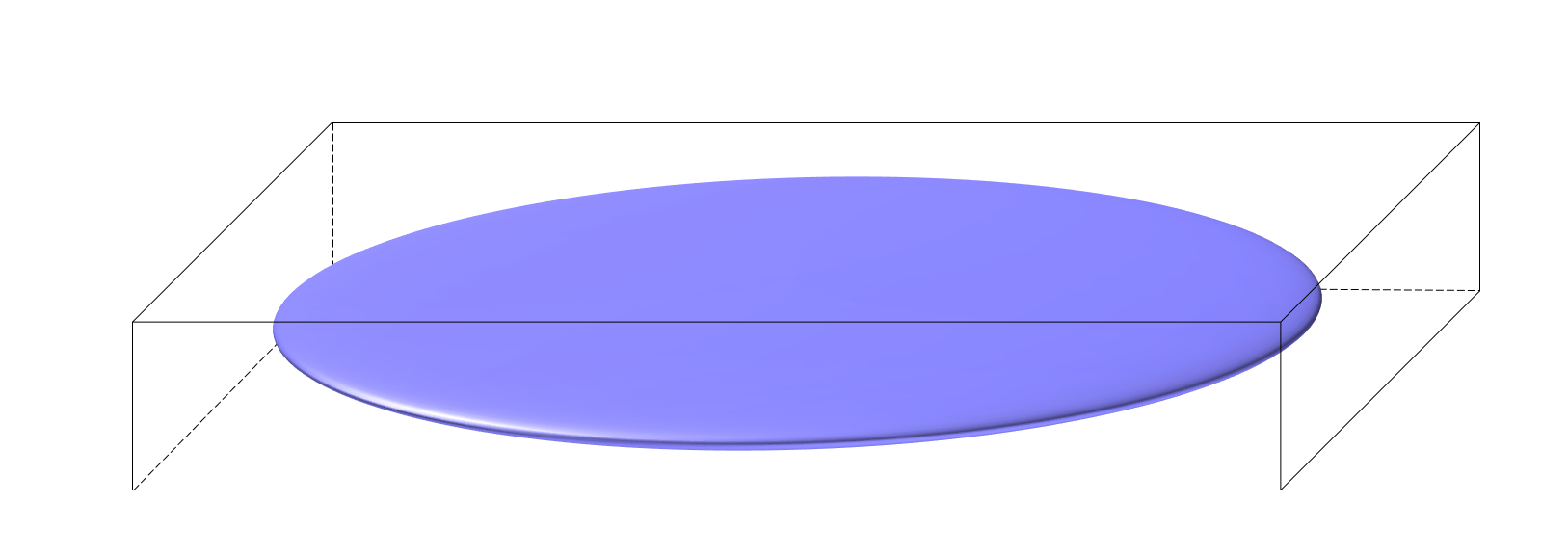}
\caption{Cigar-shaped (left) and pancaked-shaped (right) densities. }
\label{Anisotropy}
\end{figure}

\vspace{0.3em}
There are three major challenges for this convolution-type potential evaluation.
   \begin{enumerate}
    \item \textbf{Singularity:} The kernel and its Fourier transform are both singular, and sometimes the singularity of $\widehat U(\bk)$ are even stronger.
     \item \textbf{Nonlocality:} The potential value at a fixed \textsl{target} point $\bx$ depends on the density and kernel at every \textsl{source} point $\by$.
     \item \textbf{Anisotropy:}
     The density may be highly anisotropic. For example, in lower-dimensional confined quantum systems \cite{Bao2013Mathematical},
     its compact support extends much shorter in one or two directions, as shown in Figure \ref{anisoDomain} with the ``cigar-shaped'' and
     the ``pancake-shaped'' densities.
  \end{enumerate}

Since the density decays rapidly enough, it is reasonable to assume that the density is  compactly supported  (to a
fixed precision) in a bounded domain.
In computational practice, we truncate the whole space to a rectangular domain $\mathcal{D}:= \prod_{j=1}^d [-L_j,L_j]$ and discretize it with equally spaced points in each direction.
The smooth density $\rho(\bx)$ is well approximated by Fourier spectral method with spectral accuracy, achieving nearly optimal efficiency
thanks to implementation using discrete Fast Fourier Transform (FFT) \cite{Shen2011Spectral}.

The numerical problem is to compute the convolution $\Phi$ on uniform grid from discrete density values given on the same grid. The scheme is expected to achieve spectral accuracy with great efficiency, and capable of dealing with highly anisotropic density case.

Over the past decade, several fast and accurate methods have been developed based on Fourier spectral method.
The \textit{NonUniform Fast Fourier Transform} (NUFFT)-based method \cite{Jiang2014Fast, Bao2015Computing},
the first accurate fast algorithm, was implemented via the NUFFT algorithm.
Later, the \textit{Gaussian-Summation method} (GauSum) \cite{Exl2016Accurate}, emerged as the first purely FFT-based algorithm,
dealt with the singularity using a summation-of-Gaussian (SOG) approximation of the kernel that is done away from the origin.  In 2016, a much simpler method,
the \textit{kernel truncation method} (KTM), is proposed by Vico \textsl{et al} \cite{Vico2016Fast}, where they remove the singularity by truncating the kernel within a radially symmetric domain.
In KTM, one has to zero-pad the density by a factor, that is no small than $\sqrt{d}\!+\!1$~\cite{Liu2022optimal},
in each spatial direction so as to capture the unpleasant oscillations brought by the discontinuous kernel truncation.
Unfortunately, both memory requirement and computation costs scale linearly with the anisotropy strength, and it places a huge burden in higher space dimension, especially the three-dimension problem.

Following the kernel truncation line, Greengard \textsl{et al} \cite{Greengard2018anisotropic} proposed the \textit{anisotropic truncated kernel method} (ATKM),
where the kernel was effectively truncated within an anisotropic double-sized rectangular geometry. The computation of the Fourier transform, defined on the anisotropic domain, is even more challenging and successfully precomputed via the \textbf{near-field} approximation SOG that is performed much closer to the origin compared with GauSum. Even though GauSum and ATKM achieve nearly optimal efficiency with minimum memory requirement for both isotropic and anisotropic cases, they are technically challenging in implementation. Moreover, in practice, there are cases where high-accuracy computation is of great importance or even unavoidable, for example, the fine structure of vortices and the dynamics in context of the Bose-Einstein condensates \cite{Bao2013Mathematical}, therefore, it renders the necessity of new computational methodology.

\

It is worthy to point out that once the kernel is approximated by $U^\varepsilon(\bx)$ at the far field with great accuracy,
the Fourier transform of residual kernel $(U-U^\varepsilon)$ is quite close to
its counterpart that is defined on bounded domain. It is more favorable if the \textbf{far-field} approximation $U^\varepsilon$ is smooth.
As is noticed in classical Ewald summation method~\cite{Greengard2023dmk}, the Coulomb kernel is split into far-field and local parts as follows
\bea
U(r) := \frac{1}{r} = \frac{1}{r} \operatorname{Erf}\left(\frac{r}{\varepsilon}\right)
+ \frac{1}{r} \operatorname{Erfc}\left(\frac{r}{\varepsilon}\right)
:= U^\varepsilon(r) + (U-U^\varepsilon)(r),
\eea
where $\operatorname{Erf}(x)=\frac{2}{\sqrt{\pi}} \int_0^x e^{-t^2} {\rm d} t$ and $\operatorname{Erfc}(x) = 1\!-\!\operatorname{Erf}(x)$ are
the error and complementary error functions
and $\varepsilon$ is a positive parameter to be chosen later.
The far-field part $U^\varepsilon$ is smooth and accurately represents the kernel $1/r$ at far field,
because $\operatorname{Erf} (r/\varepsilon) \approx 1$ with more than 16 digits accurate when $r \ge 6  \varepsilon$.
The remaining near-field part $(U-U^\varepsilon)$, though singular, is compactly supported.
We split the potential into a regular integral $\Phi^{\rm{R}}(\bx)$ and a singular integral $\Phi^{\rm{S}}(\bx)$ as follows
\bea
\Phi(\bx) &=& \int_{\mathbb{R}^d} U^{\varepsilon}(\bx -\by) \rho(\by) {\rm d}\by + \int_{\mathbb{R}^d} (U-U^{\varepsilon})(\bx-\by) \rho(\by) {\rm d}\by, \notag \\
&:=& \Phi^{\rm{R}}(\bx) + \Phi^{\rm{S}}(\bx).
\label{ConvIntSplit}
\eea
Both integrals are well-defined simply because integrand of $\Phi^{\rm R}(\bx)$ is smooth and $(U- U^{\varepsilon}) \in L^1(\mathbb{R}^d)$.

The regular integral is well resolved by trapezoidal rule and implemented with nearly optimal efficiency thanks to FFT.
We choose to integrate the singular integral potential $\Phi^{\rm{S}}(\bx)$ by switching to Fourier space following
exactly the same way of ATKM~\cite{Greengard2018anisotropic}. Fortunately,
the required bounded domain Fourier transform of $(U - U^{\varepsilon})$ can be replaced by
the whole space Fourier transform, which admits analytical explicit formula, with high accuracy,
thus waiving the use of technically complicated quadrature.
Such methodology can be easily extended to anisotropic density case,
and it achieves an anisotropic strength-independent memory requirement and computation complexity.
Moreover, the numerical quadrature of each integral can be reduced to discrete convolutions, therefore,
our method is finally simplified into one discrete convolution by combining them two.
The tensor generation procedure is very simple, highly efficient (involving only FFT) and independent of the anisotropy strength.

Apart from the Coulomb kernel shown above, our method is applicable to a large class of radially symmetric kernels as long as they grow no faster than exponential functions, that is,
\beas
\int_{R_0}^{\infty} |U(r)|~  r^{d-2}~ e^{-r^2/\delta^2} {\rm d}r < \infty, \quad \delta >0,
\eeas
where $R_0 :=\min\limits_{j=1, \cdots, d} \{ 2 L_j \}$.
The far-field smooth approximation (\textbf{FSA}) function $U^\varepsilon$  only needs to satisfy the following two properties
\begin{equation}
\label{ApproxProp}~
\scalebox{0.96}{\boxed{
\begin{aligned}
&(a)~U^{\varepsilon}(\bx)~ \mbox{ is } \textit{smooth} \mbox{ and } \textit{radially symmetric}. \\[0.3em]
&(b)~\int_{R_0}^{\infty}|(U - U^{\varepsilon})(r)|  r^{d-1}   {\rm d}r \le \varepsilon_{\rm{tol}},
~\mbox{where}~ \varepsilon_{\rm{tol}}\ll 1~ \mbox{is the desired resolution}.
\end{aligned}
}}
\end{equation}
Such approximation can be realized with help of some special functions, including the error and smooth window functions \cite{SmoothWindow},
and we shall provide details later in Section \ref{KernelApprox}.
	
The paper is organized as follows. In Section \ref{NumMet}, we present a detailed description of the regular and singular integral evaluation, the discrete convolution structure and its extension to anisotropic case. Error estimates are given in Section \ref{ErrEst}. In Section \ref{KernelApprox}, we propose a simple way to construct far-field smooth approximation. Extensive numerical results are shown in Section \ref{NumericalResults} to illustrate the performance in terms of accuracy and efficiency. Some conclusions are drawn in Section \ref{sec-conclusions}.

\section{Numerical method}
\setcounter{equation}{0}
\label{NumMet}
In this section, we first focus on computation of the nonlocal potential generated by isotropic density, and discuss the anisotropic case later.
For simplicity, we assume that the density function $\rho$
is compactly supported in a square domain $\bR_L = [-L,L]^d$, which is discretized with $N$ equally spaced grid points in each spatial direction.
We fix $h = 2L/N$ as the mesh size, and
the uniform mesh grid set is denoted
\begin{equation}
\mathcal T := \left\{ (x_{1}, x_{2},\cdots, x_{d}) \Big|~ x_{j} \in \left\{h\ell,~\ell=-N/2,\cdots, N/2-1\right\}, ~j = 1, \cdots, d\right\}.
\label{UniformMesh}
\end{equation}
Following the previous discussion, we shall illustrate the computation of $\Phi^{\mathrm{R}}(\bx)$ and $\Phi^{\mathrm{S}}(\bx)$, and present a detailed analysis in terms of accuracy and efficiency.

\subsection{Evaluation of the regular integral $\Phi^{\rm{R}}(\bx)$}
Due to the compact support of the density, we have
\bea\label{RegInt-Exp}
\Phi^{\rm{R}}(\bx) = \int_{\mathbb{R}^d} U^{\varepsilon}(\bx -\by) \rho(\by) {\rm d} \by
   \approx \int_{\bR_L} U^{\varepsilon}(\bx -\by) \rho(\by) {\rm d}\by, \quad \bx \in \bR_{L}.
\eea
The above integral is well approximated by applying the trapezoidal rule quadrature \cite{Shen2011Spectral}, and
the resulting summation is reduced to a discrete convolution of tensor and the
discrete density.
The discrete convolution structure can be efficiently accelerated using FFT within ${\mathcal O}( d (2N)^d \log(2N) )$ float operations \cite{Liu2022optimal}.
For simplicity, we only present the detailed scheme for 2D case and extension to 3D case is  straightforward.
To be exact, we obtain
\bea
\Phi^{\rm R}(x_n, y_m)
&\approx&  h^2 \sum_{(n', m') \in \mathcal{I}_{N}} U^{\varepsilon}\big( (n-n') h, ~ (m-m') h \big)~ \rho(x_{n'}, y_{m'})  \notag \\
& := & \sum_{(n', m') \in \mathcal{I}_{N}} T^{(1)}_{n-n', m-m'} ~ \rho_{n', m'}, \label{Tensor1}
\eea
where $\rho_{n', m'} : = \rho(x_{n'}, y_{m'})$ and the tensor $T^{(1)}_{n, m}$ is given explicitly as
\beas
T^{(1)}_{n, m} = h^2 ~ U^{\varepsilon}(n h ,m h).
\eeas

\subsection{Evaluation of the singular integral $\Phi^{\rm S}(\bx)$}
To compute $\Phi^{\rm S}(\bx)$ within $\textbf{R}_L$,
we first reformulate it as
\bea \nonumber
\Phi^{\rm S}(\bx) &:=& \int_{\mathbb R^d} (U-U^\varepsilon)(\bx-\by)\rho(\by) {\rm d}\by \approx \int_{\textbf{R}_L} (U-U^\varepsilon)(\bx-\by)\rho(\by) {\rm d}\by \\
&=& \int_{\bx + \textbf{R}_L} (U-U^\varepsilon)(\by)\rho(\bx-\by) {\rm d}\by
= \int_{ \textbf{R}_{2L}} (U-U^\varepsilon)(\by)\rho(\bx-\by) {\rm d}\by. \label{atkm_I2}
\eea
The last identity holds true because the density is compactly supported in $\textbf{R}_L$. To be specific,
for any $\bx \in \textbf{R}_L$, we have
\begin{equation*}
   \by \in \bR_{2L}\setminus \left(\bx + \bR_L\right)
   ~~\Longrightarrow~~ \bx-\by \notin \bR_{L} ~~\Longrightarrow~~\rho(\bx-\by) = 0.
\end{equation*}
To integrate Eqn.~\eqref{atkm_I2}, one needs to obtain a good approximation of $\rho(\bx)$ over $\textbf{R}_{3L}$.
It is natural to first zero-pad the density from $\bR_L$ to $\textbf{R}_{3L}$ and then construct a Fourier series approximation therein.
Thanks to periodicity of Fourier basis, a two-fold zero-padding to $\textbf{R}_{2 L}$ is sufficient to guarantee spectral accuracy \cite{Exl2016Accurate,Liu2022optimal}, and it reads as follows
\be \label{DisFourierSeries}
\rho_{N}(\bz) := \sum_{\bk\in \Lambda} \widetilde{\rho}_{\bk} ~  e^{i \bk \cdot \bz}, \quad \bz \in \bR_{2L},
\ee
where $\Lambda =\{ \bk:=\frac{ \pi}{2 L }(k_1,\cdots, k_d) \in \frac{ \pi}{2L} \mathbb Z^d~\big| k_j= -N,\cdots, N-1, ~ j = 1,\cdots,d\}$
denotes the Fourier mesh grid.
The Fourier coefficients are defined as
\bea\label{rhoTildeK}
\widetilde \rho_\bk  = \frac{1}{(2 N)^d} \sum_{\bz_\bp \in \mathcal T}
\rho(\bz_\bp) e^{-i\bz_\bp\cdot \bk}, \quad ~\bk \in \Lambda.
\eea
As is shown earlier, the periodic extension of $\rho_N(\bz)$ is also a spectral approximation over $\textbf{R}_{3L}$, therefore, after
substituting $\rho_N$ for $\rho$ in \eqref{atkm_I2}, we obtain
\bea
\Phi^{\rm{S}}(\bx)&\approx&\int_{\textbf{R}_{2L}} (U-U^{\varepsilon}) (\by) \rho_{N}(\bx- \by) {\rm d} \by, \nonumber \\
&=&  \sum_{\bk\in \Lambda} ~\left[\int_{\textbf{R}_{2L}} (U-U^{\varepsilon})(\by) e^{-i \bk \cdot \by} { \rm d} \by \right]~ \widetilde{\rho}_{\bk} ~e^{i \bk \cdot \bx}, \notag \\
\quad\quad &:=& \sum_{\bk\in \Lambda} W(\bk) ~\widetilde{\rho}_{\bk} ~e^{i \bk \cdot \bx}, \quad \bx \in \bR_L,\label{I_num2}
\eea
where $W(\bk)$ denotes the Fourier transform of $(U-U^\varepsilon)$ over $\bR_{2L}$ and is given explicitly
\bea
\label{Wk}
    W(\bk) = \int_{\textbf{R}_{2L}} (U-U^{\varepsilon})(\by) e^{-i \bk \cdot \by} { \rm d} \by.
\eea

Usually, the above Fourier transform does not admit explicit analytical expressions,
therefore, one has to design appropriate numerical quadrature. For example, Greengard et al.
\cite{Greengard2018anisotropic} proposes a fast and accurate quadrature by utilizing the
Gaussian-sum approximation.
In fact, we can approximate $W(\bk)$ by the Fourier transform of $(U-U^{\varepsilon})$ due to the second property of Eqn.~\eqref{ApproxProp}.
That is,
\bea
W(\bk) \approx \int_{\mathbb{R}^d}~ (U-U^{\varepsilon})(\by) ~e^{-i \bk \cdot \by} { \rm d } \by.
\label{FouTran}
\eea
With a suitable parameter $\varepsilon$, we can control the approximation error as small as any prescribed precision.
The above whole space integral is well-defined since $(U-U^\varepsilon)\in L^1(\mathbb R^d)$.
The Fourier transform of $(U-U^\varepsilon)$ is radially symmetric \cite{Vico2016Fast} and given explicitly below
\bea
\widehat{(U-U^{\varepsilon})}(\bk)=
\left\{\begin{array}{ll}
2 \pi \int_{0}^{\infty} (U-U^{\varepsilon})(r)~\operatorname{J_0}(k r) ~ r ~{\rm d}r, & d = 2,\\[0.8em]
4 \pi \int_{0}^{\infty} (U-U^{\varepsilon})(r)~\frac{\sin{(k r)}}{k r}~ r^2~{\rm d}r, & d = 3,
      \end{array}
\right.
\label{Fourier-Symmetric}
\eea
where $k = |\bk|$ and $\operatorname{J_0}(r)$ is the Bessel function of first kind with index $0$.
For common kernels, including Poisson, Coulomb and Biharmonic kernels, Eqn.~\eqref{Fourier-Symmetric} has analytical expressions.
For a more general kernel, one may resort to numerical integration, e.g., the Gauss-Kronrod quadrature,
so to obtain an accurate approximation of the Fourier transform.

\

\begin{remark}[Parameter choice of $\varepsilon$ for isotropic case]\label{rmk-parameter}
The parameter $\varepsilon$ is chosen to satisfy condition (b) of Eqn.~\eqref{ApproxProp}.
Roughly speaking, $\varepsilon < R_0/5.85$ $(R_0/8.65)$ yields about 16 (34) digits of accuracy for 3D Coulomb potential.
Numerical experience suggests that we can set $\varepsilon = 1$ to achieve 34 digits for $L = 8$ ($R_0 = 2L =16$). A detailed derivation is provided
 in Appendix \ref{epsChoose}.
\end{remark}

\

The quadrature is given on uniform mesh grid and can be rewritten as a discrete convolution of a tensor and the discrete density.
The tensor is actually the inverse discrete Fourier transform of vector $\{W(\bk)\}_{\bk \in \Lambda}\in \mathbb C^{(2N)^{d}}$.
To be exact, let us take the 2D case as an example.
Define the index set
\begin{equation}
\mathcal{I}_N= \left\{ (n, m) \in \mathbb{Z}^2 | -N/2 \le n,m\le N/2-1\right\}
\end{equation} and the Fourier modes $\mu_p = \frac{\pi p}{2L},~\mu_q=\frac{\pi q}{2L}$.
Plugging Eqn.~\eqref{rhoTildeK} into Eqn.~\eqref{I_num2} and switching the summation order, the $\Phi^{\rm{S}}$ on a uniform grid can be rewritten as
\bea
\Phi^{\rm{S}}(x_n, y_m) &=& \sum_{(p, q) \in \mathcal{I}_{2 N}} W(\mu_p,\mu_q)
~\widetilde{\rho}_{\bk} ~e^{\frac{2\pi i}{2 N} (p n+q m)} \notag \\
&=& \sum_{(n', m') \in \mathcal{I}_N} \bigg[ \frac{1}{ (2 N)^2} \sum_{(p, q) \in \mathcal{I}_{2 N}}
W(\mu_p,\mu_q)~  e^{\frac{2\pi i}{2 N} \big[p (n-n')+q (m-m')\big]} \bigg] \rho_{n', m'} \notag \\
 &:=&\sum_{(n', m') \in \mathcal{I}_N} T^{(2)}_{n-n', m-m'} ~\rho_{n', m'}, \label{Tensor2}
\eea
where the tensor $T^{(2)}_{n,m}$ can be computed out \textit{once for all} within $\mathcal{O}(4N^2 \log(4N^2))$ float operations using iFFT.

\subsection{Discrete tensor structure} \label{FastAlgo}
As pointed out in \cite{Vico2016Fast}, the KTM can be rewritten as
a discrete convolution of a tensor and the discrete density.
In fact, for any real-valued kernel function, the potential evaluation algorithm discretized on uniform grid can be rewritten as such a discrete convolution, including the aforementioned NUFFT \cite{Jiang2014Fast, Bao2015Computing}, ATKM \cite{Greengard2018anisotropic} and GauSum method \cite{Exl2016Accurate}.
Once the tensor $T$ has been generated, all algorithms share the same efficiency. The differences lie in the memory requirement and efficiency of the tensor generation procedure, and whether they are dependent on anisotropy strength.

In our method,  both the regular integral potential $\Phi^{\rm{R}}$  and the singular integral potential $\Phi^{\rm{S}}$ can be written as discrete convolutions (i.e., \eqref{Tensor1} and \eqref{Tensor2}).
Then it can be simplified to a single discrete convolution as follows
 \bea
 \label{DiscConv}
\Phi_{N}(x_n, y_m) := \sum_{(n', m') \in \mathcal{I}_N} T_{n-n', m-m'} ~\rho_{n', m'},
\eea
where $T_{n, m}=T^{(1)}_{n, m} + T^{(2)}_{n, m}$.

The tensor generation procedure is simple and efficient, and it involves purely iFFT.
Note that the tensor $T$ is axis-symmetric, i.e., $T_{-n, m}=T_{n,m} = T_{n,-m}$ due to the symmetry of
$U$ and  $U^{\varepsilon}$, therefore, we can further reduce the memory requirement to $1/4$ of its original size.
Once the tensor is available, the discrete convolution \eqref{DiscConv} can be implemented with only FFT/iFFT and pointwise multiplication on vectors
of length $(2N)^{d}$.
We refer the reader to \cite{Liu2022optimal,Vico2016Fast} for more details.

\subsection{Anisotropic density}
\label{Sec:anis}
We assume that the density is compactly supported in an anisotropic rectangle
$
\textbf{R}_{L}^{\bm{\gamma}}:= \prod_{j=1}^d  [-L \gamma_j, L \gamma_j ]$, %\eea
which is also the domain of interest. We define the \textit{anisotropy vector} by $\bm{\gamma}=(\gamma_1, \cdots, \gamma_d)$. The magnitudes of the $\gamma_j$ reflect the anisotropy strength along the $j$-th direction.
Without loss of generality,
we take $\gamma_1=1, 0<\gamma_j \le 1$ for $j=2, \cdots, d$ and define the anisotropy strength as $\gamma_f := \prod_{j=1}^d \gamma_j^{-1}.$
The density is sampled on a uniform mesh grid with the same number of grid points in each spatial direction (denoted by $N$) with mesh grid  $h_j = 2 L \gamma_j/N$.

We adapt it to the anisotropic density case and derive a similar discrete convolution structure.
The difference is that the rectangular domain becomes anisotropic as $\textbf{R}_{2L}^{\bm{\gamma}}$ and the parameter $R_0$
is now
\beas
 R_0 = \min\limits_{j=1, \cdots, d} \{2L \gamma_j \}.
 \eeas
 To be specific, the discrete convolution tensor  for $d=2$ reads as
\bea
\label{TensorAnis}
\qquad
T_{n, m}
= h_1 h_2 U^{\varepsilon}(n h_1 ,m h_2 )+
\frac{1}{ (2 N)^2 } \! \sum \limits_{(p, q) \in \mathcal{I}_{2 N}}
W\Big(\frac{\pi p}{2L \gamma_1}, \frac{\pi q}{2L \gamma_2}\Big) ~ e^{\frac{2\pi i}{2 N} (p n+q m)},
\eea
where $W(\bk)$, the Fourier transform of $(U-U^\varepsilon)$ over $\bR_{2L}^{\bm{\gamma}}$, is given explicitly
\bea
\label{WkAnis}
    W(\bk) = \int_{\textbf{R}_{2L}^{\bm{\gamma}}} (U-U^{\varepsilon})(\by) e^{-i \bk \cdot \by} { \rm d} \by.
\eea
 Similar to the isotropic case, $W(\bk)$ can be approximated by the Fourier transform of $(U-U^{\varepsilon})$ due to the second property of Eqn.~\eqref{ApproxProp} with a suitable $\varepsilon$.

\

\begin{remark}[Parameter choice of $\varepsilon$ for anisotropic case]
The choice of parameter $\varepsilon$ follows exactly the same principles as shown in Remark~\ref{rmk-parameter}. Clearly, the parameter $\varepsilon$
decreases as the anisotropy strength increases.
Numerical experience suggests that
we can set $\varepsilon = 0.4$ to achieve 34 digits of accuracy for 3D
Coulomb potential with anisotropy strength $\bm{\gamma}=(1,1,\gamma_3)$,
 $\gamma_3=1, 1/2,1/4,1/8$ and $L=12$.
\end{remark}

\

\begin{remark}[Comparison of the tensor generation procedure]
The memory requirement and computational costs for our method (denoted by FSA hereafter), KTM, GauSum and ATKM are
$\mathcal{O}(S N^d)$ and ${\mathcal O}(S N^d \log(S N^d))$ respectively, with the only differences lying in the zero-padding factor $S$.
For FSA, GauSum and ATKM, $S=2^d$ for both isotropic and anisotropic cases.
While, in KTM, we take $S=3^d$ for isotropic case and
\begin{equation}
\label{ZeroPad-KTM}
S=\prod_{j=1}^d S_j, \quad ~ S_j :=\lceil{1\!+\gamma_j^{-1}(1+ \!\prod_{k=2}^d \gamma_k^2)^{1/2} \rceil}
\end{equation} for anisotropic case \cite{Liu2022optimal}, where $\lceil \cdot \rceil$ is the rounding up function.
\end{remark}

\section{Error estimates}
\label{ErrEst}
In this subsection, we choose to present the error estimates, which mainly consist of errors coming from trapezoidal rule discretization,
the density approximation, and Fourier transform of the kernel on rectangular domain. For more details, we refer the readers to
Appendix \ref{App-ErrEst}.

\

\begin{thm} \label{FourierApproximation}
For smooth density $\rho(\mathbf{x})$ that is compactly supported in $\bm{\mathrm{R}}_L^{\bm{\gamma}}$, the following estimates
  \beas
  \|\Phi -\Phi_{N} \|_{L^{\infty}(\bm{\mathrm{R}}_{L}^{\bm{\gamma}})}  &\lesssim &  N^{-m} + \varepsilon_{\rm{tol}},
  \eeas
hold true for any positive integer $m >0$.
\end{thm}
\section{Far-field smooth approximation} \label{KernelApprox} 
\begin{figure}[!htbp]
  \label{smothkernel}
\centering
\includegraphics[width=6.4cm,height=4.8cm]{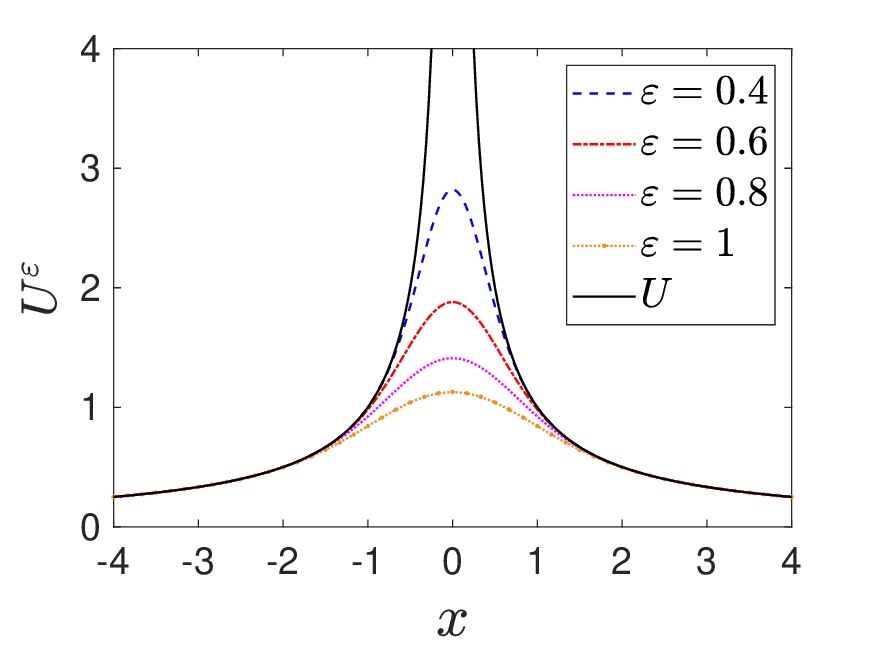}
\includegraphics[width=6.4cm,height=4.8cm]{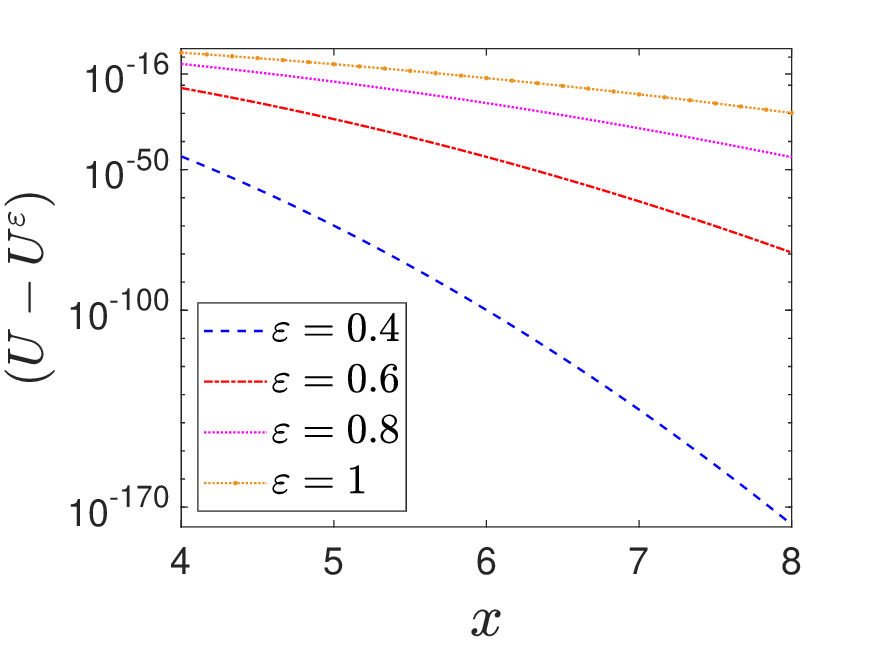}
\caption{Slice plot of $U^{\varepsilon}$ (left) and the error $(U-U^{\varepsilon})$ (right) for Coulomb kernel $U = 1/r$. }
\end{figure}

\setcounter{equation}{0}
Notice that the smooth approximation is done far away from the origin, i.e.,  $[R_0, \infty]$, and $\operatorname{Erf}(r/\varepsilon) \approx 1$ with more than 16 (34) digits of accuracy for $r \ge 6\varepsilon$ ($8.7\varepsilon$),
as shown in Figure \ref{smothkernel} for the Coulomb kernel, we can define the approximation as
$
U^{\varepsilon}(r):= U(r) \operatorname{Erf}(\frac{r}{\varepsilon}).
$
 Such approximation satisfies the requirement (b) of Eqn.~\eqref{ApproxProp} and a rigorous
 derivation is provided in Theorem \ref{Theorem:Approx} below.

\begin{thm}[\textbf{Far-field approximation}]
\label{Theorem:Approx}
For $d$-dimensional radially symmetric kernel $U(r)$ that satisfies
$ \int_{R_0}^{\infty} |U(r)|~  r^{d-2}~ e^{-r^2/\delta^2} {\rm d}r < \infty$ ($\delta>0$) and its approximation
\bea\label{KerAppr}
U^{\varepsilon}(r)= U(r) \operatorname{Erf}\left(\frac{r}{\varepsilon}\right),
\eea
there exists a positive constant $\varepsilon_0>0$ such that following estimate
\bea
\label{ApproInteg}
\int_{R_0}^{\infty}|(U - U^{\varepsilon})(r)|~ r^{d-1}~ {\rm d} r \le \varepsilon_{\rm{tol}},\quad \forall ~0 < \varepsilon < \varepsilon_0
\eea
holds true.
\end{thm}

\begin{myproof}
Using asymptotics of complementary error function, we obtain
\beas
 \operatorname{Erfc}(x) &=& \frac{2}{\sqrt{\pi}} \int_x^{\infty} e^{-t^2} {\rm d}t = \frac{2}{\sqrt{\pi}} \int_x^{\infty} e^{-t^2+x^2-x^2} {\rm d}t = \frac{2}{\sqrt{\pi}} e^{-x^2} \int_x^{\infty} e^{-(t-x)(t+x)} {\rm d}t \\
 &\le& \frac{2}{\sqrt{\pi}} e^{-x^2} \int_x^{\infty} e^{-2 x(t-x)} {\rm d}t = \frac{1}{\sqrt{\pi}} \frac{1}{x} e^{-x^2}.
\eeas
Then, we have
\beas
\int_{R_0}^{\infty} |(U - U^{\varepsilon})(r)| r^{d-1} {\rm d} r &=&
\int_{R_0}^{\infty} |U(r)|  r^{d-1} \operatorname{Erfc}\left(\frac{r}{\varepsilon}\right) {\rm d} r
\le \frac{\varepsilon}{\sqrt{\pi}} \int_{R_0}^{\infty} |U(r)| r^{d-2} e^{-\frac{r^2}{ \varepsilon^2}} {\rm d}r.
\eeas
Function $f(\varepsilon): = \int_{R_0}^{\infty} |U(r)| ~ r^{d-2} ~ e^{-r^2/\varepsilon^2}  ~{\rm d} r$ is monotone increasing and bounded.
Using dominated convergence theorem, we obtain
\beas
\lim_{\varepsilon \rightarrow 0} f(\varepsilon) = \lim_{k \rightarrow \infty} f(\varepsilon_k)
:=\lim_{k \rightarrow \infty} f(\delta/k)  = \int_{R_0}^{\infty} |U(r)| ~ r^{d-2} \lim_{k \rightarrow \infty} e^{-\frac{r^2 k^2}{\delta^2}}  {\rm d} r = 0,
\eeas
 Therefore, there exists $\varepsilon_0>0$ such that estimate \eqref{ApproInteg} holds true for any $\varepsilon < \varepsilon_0$.
\end{myproof}

\vspace{0.5em}

As for the \textbf{smoothness} requirement, i.e., condition (a) of Eqn.~\eqref{ApproxProp},
the approximation \eqref{KerAppr} is not necessarily smooth for all kernels.
For non-smooth cases, we can first set $U^{\varepsilon}(r)$ to be $0$ for $r\approx 0$ and $U(r) \operatorname{Erf}(r/\varepsilon)$ for $r \ge R_0$,
then connect the two regions smoothly via a smooth vanishing window function \cite{SmoothWindow}.  To be specific,
we set
\bea
\label{GeneralApprox1}
U^{\varepsilon}(r) = U(r) \operatorname{Erf}\left(\frac{r}{\varepsilon}\right) \big[1-\xi(r)\big],
\eea
where $\xi(r) = 0$ if formula \eqref{KerAppr} is smooth, otherwise,
set $\xi(r) $ to be a smooth vanishing window function.
Such construction mechanism works for almost all kernels, but it may not be the simplest one for certain kernels.
For example, for the 2D Poisson kernel $U(r) = -\frac{1}{2\pi} \ln(r)$, we choose the following approximation function,
\bea
\label{lnrAppr}
U^{\varepsilon}(r) = -\frac{1}{2 \pi} \Big[\ln(r) + \frac{1}{2} \operatorname{E_1}
\Big(\frac{r^2}{\varepsilon^2}\Big) \Big],
\eea
where $\operatorname{E_1}(r):= \int_r^{\infty} t^{-1} e^{-t} {\rm d}t $  is the exponential integral function \cite{Abramowitz1964Handbook}.

\

\
With the above-mentioned smooth approximations, the Fourier transform of $(U - U^{\varepsilon})$ is reduced to one-dimensional integral \eqref{Fourier-Symmetric}.
Here we summarize the corresponding analytical expressions for different common kernels in Table \ref{tab:smoothTruFour},
from which we can observe clearly that the corresponding Fourier transform is smooth and non-oscillatory, as against the annoying oscillatory Fourier transform in KTM.

\begin{table}[htpb!]
\centering
\caption{Far-field smooth approximation and the corresponding Fourier transform for common kernels.}
\label{tab:smoothTruFour}
\begin{adjustbox}{width=1.00\textwidth}
\begin{tabular}{llll}
\toprule
  & $U(r)$ &  $U^{\varepsilon}(r)$     &$(\widehat{U - U^{\varepsilon}})(k)$\\
\midrule
\rule{0pt}{16pt}
\footnotesize{\textbf{Poisson}} &  {\footnotesize{$2D$:}} $ \frac{-1}{2 \pi} \ln(r) $ &  $\frac{-1}{2 \pi} \big[\ln(r)+\frac{1}{2} \operatorname{E_1}(\frac{r^2}{\varepsilon^2}) \big] $ & $ \frac{1}{k^2} \big[1-e^{-\frac{1}{4} k^2 \varepsilon^2}\big]$ \\[6pt]
\midrule
\rule{0pt}{16pt}
\multirow{2}{*}{\footnotesize{\textbf{Coulomb}}}
 &  {\footnotesize{$2D$:}} $ \frac{1}{2 \pi r} $ & $ \frac{1}{2 \pi r} \operatorname{Erf}(\frac{r}{\varepsilon}) $ & $  \frac{1}{k}\operatorname{Erf}(\frac{k \varepsilon}{2})$ \\[6pt]
 &  {\footnotesize{$3D$:}} $\frac{1}{4 \pi r} $ & $ \frac{1}{4 \pi r} \operatorname{Erf}(\frac{r}{\varepsilon}) $ & $ \frac{1}{k^2} \big[1-e^{-\frac{1}{4} k^2 \varepsilon^2}\big]$ \\[6pt]
 \midrule
 \rule{0pt}{16pt}
 \multirow{2}{*}{\footnotesize{\textbf{Biharmonic}}}
 & {\footnotesize{$2D$:}} $ \frac{- r^2 }{8 \pi  } \big[\ln(r)-1\big] $
 &  $\frac{ -r^2}{8 \pi} \big[\ln(r)+\frac{1}{2} \operatorname{E_1}(\frac{r^2}{\varepsilon^2})-1\big]$ & $  \frac{-16+e^{-\frac{1}{4} k^2 \varepsilon^2} (16+4 k^2 \varepsilon^2+k^4 \varepsilon^4)}{16 k^4}$ \\[6pt]
  & {\footnotesize{$3D$:}} $ \frac{r}{8 \pi}$
 & $\frac{r}{8 \pi} \operatorname{Erf}(\frac{r}{\varepsilon}) $ & $  \frac{-8+e^{-\frac{1}{4} k^2 \varepsilon^2} (8+ 2 k^2 \varepsilon^2+k^4 \varepsilon^4)}{8 k^4}$ \\[6pt]
\bottomrule
\end{tabular}
 \end{adjustbox}
\end{table}
The parameter $\varepsilon$ is chosen to satisfy condition (b) of Eqn.~\eqref{ApproxProp}, ant it depends on kernel $U$, far-field approximation $U^{\varepsilon}$,
$R_0$ and $\varepsilon_{\rm{tol}}$.
A detailed derivation is provided in Appendix \ref{epsChoose}, using which we can
set $\varepsilon = 1 (0.4)$ for the isotropic (anisotropic) density case for the 2D/3D Coulomb and 2D Poisson potentials when $\varepsilon_{\rm{tol}} =  10^{-16}$ and $L = 12$.

\begin{remark}
In computational practice, we shall choose a relatively large $\varepsilon$ for better numerical accuracy. Even though the integrand of regular integral \eqref{RegInt-Exp} is smooth for any fixed $\varepsilon$ on the continuous level,
for a smaller $\varepsilon$, it requires a finer mesh so to capture the sharp variations near the singularity on the discrete level.
\end{remark}

\section{Numerical results}
\label{NumericalResults}
\setcounter{equation}{0}

In this section, we shall investigate the accuracy (in both double and quadruple precision) and efficiency for different nonlocal potentials evaluation with isotropic and anisotropic densities.
The computational domain $\bR_{L}^{\bm \gamma}$ is discretized uniformly in each spatial direction with mesh size $h_j$,
and we define  mesh size vector as $\bh = (h_1, \cdots, h_d)$.
For simplicity, we shall use $h$ to denote the mesh size if all the mesh sizes are equal.

The numerical error is measured in following norm
\begin{equation*}
{\mathcal E}:=\|\Phi-\Phi_{\bh}\|_{l^\infty}/\|\Phi\|_{ l^\infty}
=\textrm{max}_{\bx \in \mathcal T_{\bh}} |\Phi(\bx)-\Phi_{\bh}(\bx)|/\textrm{max}_{\bx \in \mathcal T_{\bh}} |\Phi(\bx)| ,
\end{equation*}
where $\Phi_{\bh}$ is the numerical solution on mesh grid $\mathcal T_{\bh}$ and $\Phi(\bx)$ is the exact solution.
In the following examples, the potential $\Phi(\bx)$ can be computed analytically. For simplicity, we do not list them here but refer to \cite{Greengard2018anisotropic, Liu2022optimal}.

\

\subsection{The Coulomb potentials in 2D/3D }
\begin{exmp} \label{3DPoissonInteraction}
Here, we consider the 2D/3D Coulomb potentials with the following two types of source density
\begin{itemize}
\item {\bf Case I}:  Isotropic/anisotropic density
\beas
\rho(\mathbf{x})=
\left\{\begin{array}{ll}
e^{-(x^2+y^2/\gamma_2^2)/\sigma^2},  &d=2,  \\[0.5em]
e^{-(x^2+y^2+z^2/\gamma_3^2)/\sigma^2}, & d=3. \\
\end{array}\right.
\eeas
\item {\bf Case II}:  Shifted  density  $\rho(\mathbf{x})=\rho_0(\mathbf{x})+ \rho_0(\mathbf{x}-\mathbf{x}_0)$, where $\rho_0(\mathbf{x})$ is generated by taking the Laplacian of the potential $\Phi_0(\mathbf{x}) = e^{-(x^2+y^2+z^2/ \gamma_3^2)/\sigma^2}$, i.e.,
    $
    \rho_0(\mathbf{x}) = -\Delta \Phi_0(\mathbf{x}) = \Phi_0(\mathbf{x}) \Big( -\frac{4x^2}{\sigma^4}- \frac{4y^2}{\sigma^4}-\frac{4z^2}{ \gamma_3^4 \sigma^4} + \frac{4}{\sigma^2}+ \frac{2}{\gamma_3^2 \sigma^2} \Big).
    $
\end{itemize}
\end{exmp}

\

Table \ref{Coulomb3Dtkm} presents errors of the 2D/3D Coulomb potentials in double and quadruple precision for isotropic and anisotropic densities.
For isotropic case, i.e., $\gamma_2 = 1$ or $\gamma_3 = 1$, we consider Case I with $L = 8,~ \sigma = \sqrt{0.8}$
and $\varepsilon =1$.
For anisotropic case, we consider Case I with $L=8, ~\sigma=\sqrt{1.2}, ~\varepsilon =0.5, ~\bh=\bm{\gamma}/4$
and Case II with $L=12,~ \sigma=\sqrt{0.8}, ~\varepsilon =0.4,~ \bx_0=(1,1,0)$ and $\bh=\bm{\gamma}/8$, where
$\bm{\gamma} = (1, \gamma)$ for $d=2$ and $\bm{\gamma} = (1, 1, \gamma)$ for $d=3$.

From Table \ref{Coulomb3Dtkm}, we can conclude that our method achieves spectral accuracy for both isotropic and anisotropic cases.
Note that it can reach $10^{-34}$ for quadruple-precision version when
the tensor $T$ in \eqref{DiscConv} is generated with high precision.

\begin{table}[!htbp]
\centering
\caption{Errors of the 2D/3D  Coulomb potentials for isotropic and anisotropic densities in Example \ref{3DPoissonInteraction}.}
\label{Coulomb3Dtkm}
\begin{tabular}{cccccc}
\toprule
 \multicolumn{6}{c} {Isotropic density}  \\
 \cline{2-6}\rule{0pt}{12pt}
& & $h=1$ & $ h=1/2 $ &$h= 1/4 $ & $ h=1/8 $ \\
\midrule
2D& double &  1.3856E-02 & 2.9648E-08 &  2.8012E-16 & 5.6025E-16 \\
\midrule
\multirow{2}{*}{3D}
 & double      &  2.0681E-02 & 2.5036E-06 & 5.5511E-16  & 6.9389E-16 \\
 & quad &  2.0681E-02 & 2.5036E-06 & 4.8161E-18  & 2.4195E-34 \\
\midrule
 \multicolumn{6}{c} {Anisotropic density}  \\
  \cline{2-6}\rule{0pt}{12pt}
& & $ \gamma = 1 $ & $\gamma = 1/2$ & $ \gamma = 1/4 $ & $ \gamma = 1/8 $ \\
\midrule
2D &Case I& 4.1758E-16&  2.5550E-15 & 1.5455E-15 &  1.8119E-15  \\
\midrule
\multirow{3}{*}{3D}
&Case I& 3.7007E-16&  5.3559E-15 & 5.1651E-15 &  3.9372E-15  \\
&Case II & 6.0077E-16&  6.0289E-16 & 8.0178E-16 &  1.2020E-15  \\
&Case II (quad)& 6.9529E-34&  6.9676E-34 & 1.5629E-33 &  2.7787E-33  \\
\bottomrule
\end{tabular}
\end{table}

\

\begin{exmp} \label{Exmp:Comp_ktm_Gausum}
\textbf{Comparison with KTM and ATKM}.
To compare the performance of FSA, KTM and ATKM, we choose densities in
 Case II of Example \ref{3DPoissonInteraction}.
The computation is split into two parts: the pre-computation part (\textbf{Precomp}): computation of the tensor $\widehat{T}_{\bk}$ for KTM and FSA or $\widehat{U}_R$
 for ATKM \cite{Greengard2018anisotropic},
and the execution part (\textbf{Execution}).
The algorithms were implemented in FORTRAN and run on a single 3.00GH Intel(R) Xeon(R) Gold 6248R CPU with a 36 MB cache in Ubuntu GNU/Linux with the Intel complier ifort.
\end{exmp}

\

\begin{figure}[!htbp]
\centering
\includegraphics[width=0.56\textwidth]{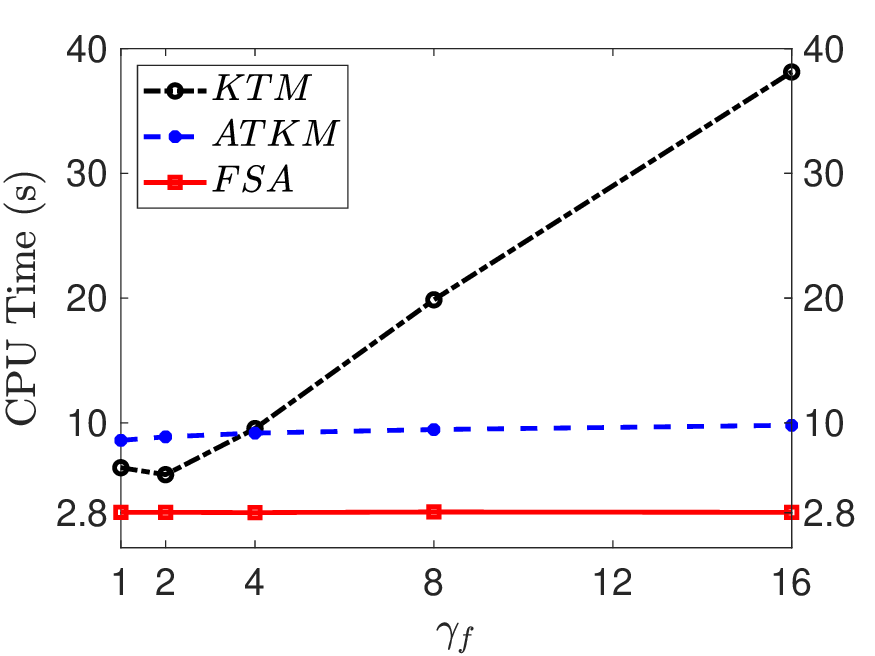}
\caption{Timing results of  the pre-computation part versus increasing anisotropy strength $\gamma_f$. }
\label{TimeCompPreCom}
\end{figure}

Table \ref{comp_ktm_SmoothTen} presents the errors $\mathcal{E}$ and CPU time for KTM and FSA in
 isotropic density case. Figure \ref{TimeCompPreCom} shows the CPU time of pre-computation part versus different anisotropy strengths $\gamma_f$.  We do not present the execution time because the subsequent computations are the same, once the pre-computation is completed.

\begin{table}[!htbp]
\centering
\caption{The performance of KTM and FSA for isotropic density.}
\label{comp_ktm_SmoothTen}
\setlength{\tabcolsep}{5mm}{
\begin{tabular}{cccc}
\toprule
$N^d= 192^3$ & $\mathcal{E}$ & $T_{\rm Precomp}(s)$ &$T_{\rm Execution}(s)$  \\
 \midrule
 KTM & 3.3502E-16 & 6.41  &2.24   \\
FSA & 3.3394E-16 &2.84   & 2.24   \\
\bottomrule
\end{tabular}}
\end{table}

From Table \ref{comp_ktm_SmoothTen} and Figure \ref{TimeCompPreCom}, we can see that:
 (i) The pre-computation of FSA is the fastest for both isotropic and anisotropic cases.
 This is because both FSA and ATKM require FFT on vectors of length $(2^d) N^d$,
 whereas KTM requires FFT on vectors of length $S N^d$ with the zero-padding factor $S$ (taken as \eqref{ZeroPad-KTM})
 growing linearly with the anisotropy strength $\gamma_f$ \cite{Liu2022optimal}.
The complicated evaluation of kernel's Fourier transform over the rectangle domain
 in ATKM results in a slightly degraded efficiency.
 (ii) The CPU time of FSA and ATKM are independent on the anisotropy strength $\gamma_f$.
 In contrast, the time for KTM exhibit a linear growth with respect to $\gamma_f$.
 Consequently, we can conclude that FSA  performs better in terms of efficiency.

\subsection{The Poisson potential in 2D}
\begin{exmp}\label{exmp:Poisson}
 We consider  the 2D Poisson potential with  density $\rho(\mathbf{x}) = e^{-|\mathbf{x}|^2/ \sigma ^2}$ in isotropic case and density
 $
\rho(\mathbf{x})= -\Delta \Phi(\mathbf{x}) =  \Phi(\mathbf{x}) \big(-\frac{4 x^2}{\sigma^4} - \frac{4 y^2}{\gamma_2^4 \sigma^4}+ \frac{2}{\sigma^2} +\frac{2}{\gamma_2^2 \sigma^2} \big)
$
 in anisotropic case,
 where $\Phi(\mathbf{x}) = e^{-(x^2+y^2/\gamma_2^2)/\sigma^2}$.
\end{exmp}

Table \ref{tab:Poisson} presents errors of the 2D Poisson potentials for isotropic and anisotropic densities.
For isotropic case, the parameters are chosen as $L=8$, $\sigma=\sqrt{1.2}$ and $\varepsilon = 1$.
  For anisotropic case, we choose $L=10,~ \sigma=1.2, ~\varepsilon = 0.4$ and $\textbf{h}=(1,\gamma_2)/8$.

  \begin{table}[!htbp]
\centering
\caption{Errors of the 2D  Poisson potentials for isotropic and anisotropic  densities in Example \ref{exmp:Poisson}.}
\label{tab:Poisson}
\begin{tabular}{ccccc}
\toprule
\multirow{2}{*}{Isotropy}
&$ h=2 $     & $h=1$ & $ h=1/2 $ &$ h=1/4 $ \\
 & 2.1786E-01&  1.3761E-03 & 5.5617E-09 &  4.9577E-16 \\
\midrule
\multirow{2}{*}{Anisotropy}
&  $ \gamma_2=1 $ & $\gamma_2=1/2$ & $ \gamma_2=1/4 $ & $ \gamma_2=1/8 $ \\
& 4.5519E-16&  2.2204E-16 & 6.2728E-16 & 1.5016E-15  \\
\bottomrule
\end{tabular}
\end{table}

\subsection{Application to some common potentials}
We test the accuracy using some common potentials, such as the 3D/quasi-2D dipole-dipole interaction (DDI), 2D/3D Biharmonic and 2D/3D Yukawa potentials.
\begin{exmp}
\label{exmp:CommonPoten}
Here, we consider the following potentials
\begin{itemize}[itemsep=8pt, topsep=4pt]
\item {\bf 3D DDI}:  The kernel is given as
\beas
 U(\mathbf{x})  =
  \frac{3}{4 \pi} \frac{\textbf{m} \cdot \textbf{n} - 3 (\mathbf{x} \cdot \textbf{m}) (\mathbf{x} \cdot \textbf{n})/|\mathbf{x}|^2 }{|\mathbf{x}|^3},
\eeas
 where $\textbf{n}$, $\textbf{m} \in \mathbb{R}^3$ are unit vectors, and
the 3D DDI is reformulated as
  $
  \Phi(\mathbf{x})= -(  \textbf{m} \cdot \textbf{n} ) \rho(\mathbf{x})-3 \frac{1}{4 \pi |\mathbf{x}|} \ast (\partial_{\textbf{n} \textbf{m}} \rho),
$
where $\partial_{\textbf{m}} = \textbf{m} \cdot  \nabla $ and $ \partial_{\textbf{n}\textbf{m}}= \partial_{\textbf{n}} (\partial_{\textbf{m}})$.
In fact, it can be calculated via the 3D Coulomb potential with source term $\partial_{\textbf{n}\textbf{m}} \rho$.
\item{\bf Quasi-2D DDI}: The kernel is given as
\begin{small}
 \beas
  U(\mathbf{x})&=&-\frac{3}{2} ( \partial_{\bn_\perp\bn_\perp}-
n_3^2\nabla_\perp^2)\frac{1}{(2\pi)^{3/2}}\int_{\mathbb R}
\frac{e^{-s^2/2}}{\sqrt{|\mathbf{x}|^2+\eta^2s^2}}{\rm d}s
:=-\frac{3}{2} ( \partial_{\bn_\perp\bn_\perp}-
n_3^2\nabla_\perp^2)  \widetilde{U}(\mathbf{x}),
\eeas\end{small}where  $\bn_\perp=(n_1,n_2)^{\top}$,
$\partial_{\bn_\perp}=\bn_\perp\cdot\nabla_\perp$ and $\partial_{\bn_\perp\bn_\perp}=\partial_{\bn_\perp}(\partial_{\bn_\perp})$. The quasi-2D DDI
can be reformulated  as
$
   \Phi(\mathbf{x}) =\widetilde U \ast \left[ -\frac{3}{2} ( \partial_{\bn_\perp\bn_\perp}-
   n_3^2\nabla_\perp^2)  \rho \right ].
$
\item {\bf Biharmonic potential}: The kernel is given as
\begin{equation*}
      U(\mathbf{x}) =
      \left\{\begin{array}{ll}
-\frac{1}{8\pi} |\mathbf{x}|^2 \left(\ln(|\mathbf{x}|)-1\right),  &d=2,  \\[0.5em]
\frac{1}{8\pi} |\mathbf{x}|, & d=3. \\
\end{array}\right.
   \end{equation*}
   \item {\bf Yukawa potential}: The kernel is given as
 \begin{equation*}
      U(\mathbf{x}) =
      \left\{\begin{array}{ll}
\frac{1}{2 \pi} \operatorname{K_0} (\lambda |\mathbf{x}|),  &d=2,  \\[0.5em]
\frac{1}{4 \pi |\mathbf{x}|}e^{-\lambda |\mathbf{x}|}, & d=3. \\
\end{array}\right.
   \end{equation*}
\end{itemize}
\end{exmp}

In our computation, we choose the density $\rho(\bx) = e^{-|\bx|^2/\sigma^2}$. Unless otherwise specified, the parameters are chosen
$L=12, ~\sigma=\sqrt{1.2}$ and $\varepsilon=1$.
Table \ref{tab:CommonPoten} shows errors of
the quasi-2D DDI, computed with $\eta=1/\sqrt{32}$ and $\bn = (0,0,1)^{\top}$,
the 3D DDI, computed with $L = 8$
with dipole orientations $\textbf{n}= (0.82778, 0.41505, -0.37751)^{\top}, \textbf{m}= (0.3118, 0.9378, -0.15214)^{\top}$,
the 2D/3D Biharmonic
and Yukawa potentials.
From Table \ref{tab:CommonPoten}, we can conclude that FSA achieves spectral accuracy.

\begin{table}
\centering
\caption{Errors of the quasi-2D/3D  DDI, 2D/3D Biharmonic and Yukawa potentials in Example \ref{exmp:CommonPoten}.}
\label{tab:CommonPoten}
\begin{tabular}{cccccc}
\toprule
 & & $ h = 2 $ & $h = 1$ & $h= 1/2 $ & $ h=1/4 $ \\
  \midrule
  \multirow{2}{*}{DDI}
 &quasi-2D & 2.0847E-01&  7.4038E-03 & 2.2647E-07 &  5.0826E-15   \\
 & 3D      & 2.2087    &  3.3668E-02 & 8.5098E-07 &  7.5667E-15   \\
\midrule
\multirow{2}{*}{Biharmonic}
 &2D  & 2.1351E-01 & 2.6558E-05 & 5.8860E-12 & 1.2938E-15 \\
 &3D  & 3.4293E-01 & 2.6307E-04 & 1.1065E-10 & 1.0623E-15 \\
 \midrule
\multirow{2}{*}{Yukawa}
&2D  & 1.7460E-01&  4.5096E-03 & 4.3501E-08 &  5.2274E-16   \\
&3D & 2.4997E-01&  6.8294E-03 &  7.3633E-08&  9.5568E-16   \\
\bottomrule
\end{tabular}
\end{table}

\section{Conclusions}\label{sec-conclusions}

Based on a far-field smooth approximation of the kernel, we presented a simple spectral fast algorithm with nearly optimal memory requirement $\mathcal{O}(2^d N^d)$
and computational cost $\mathcal{O}(2^d N^d \log(2^d N^d))$  for calculating nonlocal potentials in both isotropic and anisotropic cases.
The smoothed kernel $U^\varepsilon$ approximates the kernel at the far field with great accuracy,
and it leads to an explicit radial symmetric substitute for the Fourier transform of $(U-U^\varepsilon)$ on bounded domain.
We split the potential into two parts, and each quadrature has discrete convolution structure.
By combing both convolution structures, our method can be simplified to a single discrete convolution with explicit tensor formulation, which can be accelerated by FFT on a double-sized vector.
It is worthy to emphasize that the tensor generation is very simple, efficient and independent of the anisotropy strength.
The performance superiorities of our method were showcased with common potentials, including the Coulomb, Poisson, Biharmonic, Yukawa potentials and DDI.

\appendix
\section{Error estimates}\label{App-ErrEst}
For the regular integral $\Phi^{\mathrm{R}}(\bx)$, the error arises from the trapezoidal rule, which is spectrally accurate because the integrand is smooth and decays exponentially \cite{Shen2011Spectral}. For the singular integral $\Phi^{\mathrm{S}}(\bx)$, we have
\begin{small}
\beas
&&|E_s(\bx)| =
\Big| \int_{\bR_{2L}^{\bm{\gamma}}} (U-U^{\varepsilon})(\by) \rho(\bx-\by) {\rm d} \by -
\sum_{\bk\in \Lambda} \widehat{(U-U^{\varepsilon})}(\bk) \widetilde{\rho}_{\bk} e^{i \bk \cdot \bx} \Big| \\
&&\le \Big| \int_{\bR_{2L}^{\bm{\gamma}}} (U-U^{\varepsilon})(\by) (\rho-\rho_N)(\bx-\by) {\rm d} \by \Big| +
\Big|\sum_{\bk\in \Lambda} \Big( W-\widehat{(U-U^{\varepsilon})}\Big)(\bk) \widetilde{\rho}_{\bk} e^{i \bk \cdot \bx} \Big|
:= I_1+I_2.
\eeas
\end{small}
The error $I_1$ is of spectral accuracy because $(U-U^{\varepsilon}) \in L^{1}(\mathbb{R}^d)$
\cite{Liu2022optimal}. For $I_2$, we obtain
\beas
I_2 &\le&  \max_{\bk \in \Lambda} \Big| \int_{\mathbb{R}^d \setminus \bR_{2L}^{\bm{\gamma}} } (U-U^{\varepsilon})(\by) e^{-i \bk \cdot \by} { \rm d } \by \Big|
\sum_{\bk\in \Lambda} |\widetilde{\rho}_{\bk}|
\lesssim   \int_{\mathbb{R}^d \setminus \bR_{2L}^{\bm{\gamma}} } \left| (U-U^{\varepsilon})(\by) \right|  { \rm d } \by \\
&\le& \int_{\mathbb{R}^d \setminus {\mathbf{B}_{R_0}} }  \left| (U-U^{\varepsilon})(\by) \right| { \rm d } \by =  S^{d-1} \int_{R_0}^{\infty} \left| (U-U^{\varepsilon})(r) \right| r^{d-1} {\rm d}r \lesssim  \varepsilon_{\rm{tol}},
\eeas
where  ${\mathbf{B}_{R_0}}$ is a ball centered at the origin with radius $R_0$ and $S^{d-1} = 2 \pi^{\frac{d}{2}}/\Gamma(d/2)$.
Therefore, our method can achieve spectral accuracy with a suitable $\varepsilon$.

\section{The choice of parameter $\varepsilon$}
\label{epsChoose}
 The parameter $\varepsilon$ is chosen to guarantee condition (b) of Eqn.~\eqref{ApproxProp} is satisfied, and it varies with the kernel.
For the 3D Coulomb kernel $U(r) =1/(4\pi r)$, we have
\beas
&&\int_{R_0}^{\infty} |(U - U^{\varepsilon})(r)| ~ r^2~ {\rm d} r =
\frac{1}{4\pi} \int_{R_0}^{\infty} r \operatorname{Erfc}\left(\frac{r}{\varepsilon}\right) {\rm d} r  =  \frac{1}{4\pi}\varepsilon^2  \int_{R_0/\varepsilon }^{\infty} s \operatorname{Erfc}\left(s\right) {\rm d}s \\
&&=R_0^2 \Big[  \frac{1}{4\pi c^2}   \int_{c}^{\infty} s \operatorname{Erfc}\left(s\right) {\rm d}s \Big] : = R_0^2 f(c), \quad  \text{where} ~ c:= R_0/\varepsilon.
\eeas
We choose $ c $ such that $ f(c) R_0^2 \le \varepsilon_{\rm{tol}}$,  i.e., $f(c) \le \varepsilon_{\rm{tol}}/R_0^2$.
Clearly, $f(c)$ is monotonically decreasing and the critical value of $c$ depends on $R_0$ and  $\varepsilon_{\rm{tol}}$. For example, for $R_0 = 24$,
we can set $c \ge 5.85 $ ($8.65$), i.e.,
$\varepsilon \le R_0/5.85$ ($R_0/8.65$) to achieve 16 (34) digits of accuracy.
Similarly, for the 2D Coulomb and Poisson potentials, we can set $\varepsilon \le R_0/5.64$ and
$\varepsilon \le R_0/5.75$  respectively to guarantee 16 digits of accuracy.

\end{document}